\documentclass{article}

 \pdfoutput=1

\usepackage{hyperref}

\usepackage{fullpage}
\usepackage{amssymb}
 \usepackage{amsthm}
\usepackage{amsmath}
\usepackage[figuresright]{rotating}
 \usepackage[cp1250]{inputenc}
\usepackage{amsfonts}
\usepackage{amsmath}
\usepackage{amsthm}
\usepackage{amssymb}
\usepackage{mathrsfs}
\usepackage{graphicx}   
\usepackage{enumitem}
\usepackage{relsize}
\usepackage[colorinlistoftodos]{todonotes}

\usepackage{color}

\usepackage{epstopdf}

\theoremstyle{definition}
\newtheorem{example}{Example}[section]
\newtheorem{definition}[example]{Definition}
\newtheorem{theorem}{Theorem}[section]
\newcommand{\dd}{{\rm d}}
\newcommand{\diver}{{\rm div}}
\newcommand\DT[1]{\mathchoice
{{\buildrel{\hspace*{.1em}\text{\LARGE.}}\over{#1}}}
{{\buildrel{\hspace*{.1em}\text{\Large.}}\over{#1}}}
{{\buildrel{\hspace*{.1em}\text{\large.}}\over{#1}}}
{{\buildrel{\hspace*{.1em}\text{\large.}}\over{#1}}}}

\newcommand{\R}{\mathbb{R}}
\newcommand{\N}{\mathbb{N}}
\renewcommand{\O}{\Omega}
\newcommand{\e}{\varepsilon}

\renewcommand{\a}{\alpha}

\newcommand{\f}{\psi}
\renewcommand{\d}{{\rm d}}

\newcommand{\id}{\mbox{{\rm id}}}
\newcommand{\bulet}{\mbox{\thinspace$_{^{\mbox{$^\bullet$}}}$}}

\newcommand\llangle{\mathchoice{\big\langle\hspace{-.3em}\big\langle}
{\langle\hspace{-.2em}\langle}{\langle\!\langle}{\langle\!\langle}}
\newcommand\rrangle{\mathchoice{\big\rangle\hspace{-.3em}\big\rangle}
{\rangle\hspace{-.2em}\rangle}{\rangle\!\rangle}{\rangle\!\rangle}}

\newcommand{\Gstatic}{G}

\newcommand{\Gm}{{\mathscr G}}
\newcommand{\ccoupl}{\vec{a}}
\newcommand{\mathcalI}{\Theta}
\newcommand{\A}[1]{\langle#1\rangle}

\newcommand{\rca}{{\rm rca}}

\newcommand{\be}{\begin{eqnarray}}

\newcommand{\ee}{\end{eqnarray}}

\newcommand{\reviewMK}[1]{{\color{black}#1}}

\begin{document}

\title{Computational modeling of magnetic hysteresis with thermal effects}

\author{Martin Kru\v{z}\'{\i}k
\thanks{Institute of Information Theory and Automation of the CAS, Pod vod\'{a}renskou
v\v{e}\v{z}\'{\i}~4, CZ-182~08~Praha~8, Czech Republic}, 
Jan Valdman
\thanks{Institute of Mathematics and Biomathematics, Faculty of Science, 
University of South Bohemia, Brani\v sovsk\' a 31, CZ--37005, Czech Republic 
and Institute of Information Theory and Automation, Academy of Sciences, 
Pod vod\'{a}renskou v\v{e}\v{z}\'{\i}~4, CZ--18208~Praha~8, Czech Republic.}
}
  
\maketitle

\begin{abstract}
We study computational behavior  of a   mesoscopic model describing temperature/external magnetic field-driven evolution of magnetization.
 Due to nonconvex anisotropy energy describing magnetic properties of a body, magnetization can develop fast spatial oscillations creating  complicated microstructures. These microstructures are encoded in Young measures, their first moments then identify macroscopic magnetization. Our model assumes that changes of magnetization can contribute to dissipation and, consequently, to variations of the body temperature affecting the length of magnetization vectors. In the ferromagnetic state, minima of the anisotropic energy density  depend on temperature and they tend to zero  as we approach   the so-called Curie temperature. This   brings the specimen to a paramagnetic state. Such a  thermo-magnetic model is fully discretized and tested on two-dimensional examples. Computational results qualitatively agree with experimental observations. The own MATLAB code used in our simulations is available for download.
\end{abstract}

 {Keywords: dissipative processes, hysteresis, micromagnetics, numerical solution, Young measures }

\section{Introduction}
\label{Introduction}
In the isothermal situation, the configuration of a rigid ferromagnetic body occupying a bounded domain $\O\subset\R^d$ is usually  described by a magnetization $m:\O\to\R^d$ which denotes  density of magnetic spins and  which vanishes if the temperature $\theta$ is above the so-called  Curie temperature $\theta_{\rm c}$. 
Brown \cite{brown3} developed a theory called ``micromagnetics''  relying on the assumption that equilibrium states of saturated ferromagnets are  minima of an energy functional. This variational theory is also capable  of    predictions  of  formation of domain microstructures. We refer  e.g.~to \cite{kruzik-prohl} for a  survey on the topic. 
Starting from a microscopic description of the magnetic energy we will continue to a mesoscopic level which is convenient for analysis of magnetic microstructures. 

On  microscopic level, the magnetic  Gibbs energy consists of
several  contributions, namely an  anisotropy energy $\int_\O\psi(m,\theta)\,\d x$, where $\psi$ is the-so called anisotropy energy density describing crystallographic properties of the material, an  exchange energy $\frac12\int_\O\e|\nabla m(x)|^2\d x$ penalizing spatial changes of the magnetization, the non-local  magnetostatic energy 
$\frac{1}{2}\int_{\R^d}\!\mu_0|\nabla u_m(x)|^2\d x$,   work done by  an external magnetic field $h$ which reads $-\int_\O h(x)\!\cdot\!m(x)\,\d x$,  and a calorimetric term $\int_\O \psi_0 \, \d x$. The anisotropic energy density depends on the material properties and defines the so-called easy axes of the material, i.e., lines along which the smallest external field is needed to magnetize fully the specimen.  There are three types of anisotropy: uniaxial, triaxial, and cubic. Furthermore, $\psi$ is supposed to be a nonnegative  function, even in its first variable, i.e., $\pm m$ are assigned the same anisotropic energy. In the magnetostatic energy, $u_m$ is the magnetostatic potential related to $m$ by the Poisson problem $\mbox{div}(\mu_0\nabla u_m{-}\chi_\O m)=0$ arising from Maxwell equations.  Here $\chi_\O:\R^d\to\{0,1\}$ denotes
the characteristic function of $\O$ and $\mu_0=4\pi\times 10^{-7} \mathrm{N}/\mathrm{A}^2$ is the permeability of vacuum.

A widely used model describing steady-state isothermal configurations is due to Landau and Lifshitz \cite{landau-lifshitz,landau-lifshitz-1} (see also e.g. Brown \cite{brown3} or Hubert and Sch\"afer \cite{hubert-schafer}), relying on minimization of Gibbs' energy  with $\theta$ as a fixed parameter, i.e.,
\begin{align}\label{micromagnetics}
\left.\begin{array}{ll}
\mbox{minimize } & \Gstatic_\e(m)
:=\displaystyle{\int_\O \Big(\psi(m,\theta)+\frac12 m{\cdot}\nabla u_m
+\frac\e2 |\nabla m|^2 - h{\cdot}m \,\d x\Big)\,\d x}\\[3mm]
\mbox{subject to } & \mbox{div}(\mu_0\nabla u_m-\chi_\O m)=0 \ \mbox{ in } \R^d \ ,\\[2mm] & m\in H^1(\O;\R^d),\ \ u_m\in H^1(\R^d),
\end{array}\right\}
\end{align}
where the anisotropy energy $\psi$ is considered in the form
\begin{equation}
\psi(m,\theta):=\phi(m)+a_0(\theta-\theta_{\rm c})|m|^2-\psi_0(\theta),
\label{def-of-psi-I}
\end{equation}
where $a_0$ determines the intensity of the thermo-magnetic coupling. To see a
paramagnetic state above Curie temperature $\theta_{\rm c}$, one should
consider $a_0>0$. The isothermal part of the anisotropy energy density $\phi:\R^d\to [0,\infty)$ typically consists of two components $\phi(m) = \phi_\mathrm{poles}(m) + b_0|m|^4$, where $\phi_\mathrm{poles}(m)$ is chosen in such a way to attain its minimum value (typically zero) precisely on lines $\{ ts_\a;\ t\in\R\}$, where each $s_\a\in\R^d$, $|s_\a|=1$ determines an  axis of easy magnetization. Typical examples are $\a=1$ for uni-axial,  $1\le\a\le 3$ for triaxial, and   $1\le\a\le 4$ for cubic magnets. We can consider a uniaxial magnet with $\phi_\mathrm{poles}(m)=\sum_{i=1}^{d-1}m_i^2$, for instance. Here,  the easy axis coincides with the $d$-th axis of the Cartesian coordinate system, i.e.,  $s_\alpha:=(0,\ldots, 1)$. On the other hand, $b_0|m|^4$ is used to ensure that, for $\theta<\theta_{\rm c}$, $\f(\cdot,\theta)$ is minimized at $ts_\a$ for $|t|^2=(\theta_{\rm c}-\theta)a_0/(2b_0)$ and that $\psi(\cdot,\theta)$ is coercive. Such energy has already been used in \cite{podio-roubicek-tomassetti}. For $\e>0$, the exchange energy $\e|\nabla m|^2$ guarantees that the problem \eqref{micromagnetics} has a solution $m_\e$. Zero-temperature limits of this model consider, in addition, that the minimizers to \eqref{micromagnetics} are constrained to be valued on the sphere with the radius $\sqrt{a_0\theta_{\rm c}/(2b_0)}$ and were investigated, e.g., by Choksi and Kohn \cite{choksi-kohn}, DeSimone \cite{desim}, James and Kinderlehrer \cite{jam-kin}, James and M\"uller \cite{jam-mu}, Pedregal \cite{pedregal0,pedregal}, Pedregal and Yan \cite{PedregalYan} and many others.

 In \cite{benesova}, the authors first consider a mesoscopic micromagnetic energy arising for setting $\varepsilon:=0$ in \eqref{micromagnetics}. Moreover, it is assumed that  changes of magnetization cause dissipation which is transformed into heat.  Increasing temperature of the specimen influences its magnetic properties. Therefore, they analyze an evolutionary anisothermal mesoscopic model of a magnetic material. The aim of this paper is to discretize this model in space and time, and to  perform numerical experiments. The plan of our work is as follows. In Section~\ref{mesoscopic} we describe the stationary mesoscopic model.  The evolutionary problem is introduced in Section~\ref{sect-evol}.  Section~\ref{numerics} provides us with a numerical approximation and some computational experiments. We finally conclude with a few remarks in Section~\ref{conclusion}. An appendix then briefly introduces an important tool for the analysis as well as for numerics, namely Young measures. 

\section{Mesoscopic description of magnetization}\label{mesoscopic}
For $\e$ small, minimizers $m_\e$ of \eqref{micromagnetics} typically exhibit fast spatial oscillations, usually called  microstructure. Indeed, the anisotropy energy, which forces magnetization vectors to be aligned with the easy axis (axes), competes with the magnetostatic energy preferring divergence-free magnetization fields. It was shown in \cite{desim} by a scaling argument that  for large domains $\O$ the exchange energy contributions becomes less and less significant in comparison with other terms and thus the so-called ''no-exchange'' formulation is a justified approximation. This generically leads, however, to  nonexistence of a minimum for uniaxial ferromagnets as shown in \cite{jam-kin} without an external field $h$. Hence, various ways to extend the notion of a solution  were developed. The idea is to capture the limiting behavior of minimizing sequences of $\Gstatic_\e(m)$ as $\e \to 0$. This leads to a ``relaxed problem'' \eqref{relaxedmicromagnetics} involving possibly  so-called Young measures $\nu$'s \cite{young} which describe fast spatial  changes of the magnetization  and can capture limit patterns.

It can be proved \cite{desim,pedregal0} that this limit configuration $(\nu,u_m)$ solves the following minimization problem involving  temperature as a parameter and  a  ``mesoscopic'' Gibbs' energy $\Gstatic$:
\begin{align}\label{relaxedmicromagnetics}
\left.\begin{array}{llr}
\mbox{minimize} &
\displaystyle{\Gstatic(\nu,m) :=\int_\O \big(\f\bulet\nu+\frac12 m{\cdot}\nabla u_m-h{\cdot}m \big)\,\d x}\ \ \ \\[2mm]
\mbox{subject to} & \displaystyle{\mbox{div} \big(\mu_0\nabla u_m-\chi_\O m \big)=0} \ \ \ \mbox{ on } \R^d,\\[2mm] & m=\id \bulet\nu\hspace{7.4em}\mbox{ on } \O,\\[2mm]
& \nu\!\in\! {\mathscr Y}^p(\O;\R^d),\ \ \ m\!\in\! L^p(\O;\R^d), \ \ \ u_m\!\in\! H^1(\R^d)\ ,
\end{array}\right\}\hspace*{-1cm}
\hspace*{-.2cm}
\end{align}
where the ``momentum'' operator ``$\,\bulet\,$'' is defined by $[\f\bulet\nu](x):=\int_{\R^d} \f(s,\theta)\nu_x(\d s)$ and similarly for  $\id:\R^d\to\R^d$ which denotes  the identity and $\nu \in {\mathscr Y}^p(\O;\R^d)$. Here, the set of  Young measures ${\mathscr Y}^p(\O;\R^d)$ can be viewed as a collection of probability measures $\nu=\{\nu_x\}_{x\in\O}$ such that $\nu_x$ is a probability measure on $\R^d$  for almost every $x\in\O$. It means that $\nu_x$ is a positive Radon  measure such that $\nu_x(\R^d)=1$. We refer to Appendix for more details on Young measures.

In \cite{benesova},  the authors  built and analyzed a mesoscopic model in anisothermal situations. A closely related thermodynamically consistent model on the microscopic level was previously introduced in \cite{podio-roubicek-tomassetti} to model a ferro/para magnetic transition. Another related microscopic model with a prescribed temperature field was investigated in \cite{banas-prohl-slodicka}. The goal of this contribution is to  discretize the model from \cite{benesova}  and test it on computational examples. In order to make our exposition reasonably self-content, we closely follow the derivation of the model presented in \cite{benesova}.   We also point out that computationally efficient numerical implementation of isothermal models  can be found in \cite{carstensen-prohl, kruzik-prohl-1,kruzik-roubicek,kruzik-roubicek-2}, where such a model was used in the isothermal variant. 

In what follows we use a standard notation for Sobolev, Lebesgue  spaces and the space of continuous functions. We denote by $C_0(\R^d)$ the space of continuous functions $\R^d\to\R$ vanishing at infinity. Further, $C_p(\R^d):=\{f\in C(\R^d);\, f/(1+|\cdot|^p)\in C_0(\R^d)\}$, and 
$C^p(\R^d):=\{f\in C(\R^d);\, |f|/(1+|\cdot|^p)\le C,\ C>0\}$.

\section{Evolution problem and  dissipation}\label{sect-evol}
If the external magnetic field $h$ varies during a time interval $[0,T]$   with a horizon $T>0$, the energy of the system and  magnetic states evolve, as well. Changes of the magnetization may cause energy dissipation. As the magnetization is the first moment of the Young measure, $\nu$, we relate the dissipation on the mesoscopic level to temporal variations  of some moments of $\nu$ and consider these moments as separate variables. This approach was already used in micromagnetics in \cite{roubicek-kruzik-1,roubicek-kruzik-2} and proved to be useful also in modeling of dissipation in shape memory materials, see e.g.~\cite{mielke-roubicek}. In view of \eqref{def-of-psi-I}, we restrict ourselves to the first two moments defining $\lambda=(\lambda_1,\lambda_2)\subset\R^d\times\R=\R^{d+1}$ giving rise to the constraint
\begin{equation}
\lambda=L\bulet\nu \ ,\qquad \text{ where }\ L(m):=(m,|m|^2)
\label{moment-constraint}
\end{equation}
and consider the specific dissipation potential depending on a ``yield set'' $S\subset\R^{d+1}$
\begin{equation}\label{form-of-zeta}
\zeta(\DT{\lambda}):= \delta^*_S(\DT{\lambda})+\frac\epsilon{q}|\DT{\lambda}|^q,\qquad q\ge2,
\end{equation}
The set $S$ determines activation threshold for the evolution of $\lambda$. It is a convex compact set containing zero in its interior.
The function $\delta_S^*\ge0$ is the Fenchel conjugate of the indicator function of $S$. Consequently, it is  convex and degree-1 positively homogeneous with $\delta_S^*(0)=0$. In fact, the first term describes purely hysteretic losses, which are rate-independent and which we consider dominant,  and the second term models rate-dependent dissipation.

In view of \eqref{def-of-psi-I}--\eqref{relaxedmicromagnetics}, the specific mesoscopic Gibbs free energy, expressed in terms of $\nu$, $\lambda$ and $\theta$, reads as
\begin{subequations}\begin{align}
g(t,\nu,\lambda,\theta):=\phi\bulet \nu+ (\theta{-}\theta_{\rm c}) \ccoupl{\cdot}\lambda-\psi_0(\theta)+\frac{1}{2}m{\cdot}\nabla u_m- h(t){\cdot}m \label{gibbs}\\
\qquad\text{with }\ m=\id \bulet \nu\quad  \end{align}
\end{subequations}
where we denoted $\ccoupl :=(0,\ldots,0,a_0)$ with $a_0$ from \eqref{def-of-psi-I} and, of course, $u_m$ again from \eqref{micromagnetics}, which makes $g$ non-local.

As done already in \cite{benesova}, we relax the constraint \eqref{moment-constraint} by augmenting the total Gibbs free energy (i.e., $\psi$ integrated over $\O$) by the term $\frac\varkappa2\|\lambda-L\bulet\nu\|^2_{H^{-1}(\O;\R^{d+1})}$ with (presumably large) $\varkappa\in\R^+$ and with $H^{-1}(\O)\cong H^1_0(\O)^*$. Thus, $\lambda$'s no longer exactly represent the ``macroscopic'' momenta of the magnetization but rather are in a position of a phase field { or an internal parameter of the model.
We define the mesoscopic Gibbs free energy $\Gm$ as 
\begin{align}\label{def-of-psi}
\hspace{-3ex}\Gm(t,\nu,\lambda, \theta):=\!\int_\O\! \Big(g(t,\nu,\lambda,\theta)+\frac\varkappa2|\nabla\Delta^{-1} (\lambda - L \bulet \nu)|^2 \Big)\,\d x 
\end{align}
with $\Delta^{-1}$ meaning the inverse of the homogeneous Dirichlet boundary-value problem for the Laplacean defined as a map  $\Delta:H^1_0(\O;\R^{d+1})\to H^{-1}(\O;\R^{d+1})$.

The value of the internal parameter may influence the magnetization of the system and vice versa and, on the other hand, dissipated energy  influences the temperature of the system, which, in turn, may affect the internal parameters. In order to capture all these effects, we employ the concept of generalized  standard materials \cite{halphen-nguyen} known from continuum mechanics and couple our micromagnetic model with the entropy balance with the rate of dissipation on the right-hand side; cf.\ \eqref{EntropyEquation}.
Then the Young measure $\nu$ is  considered to evolve quasistatically
according to the  minimization principle of the Gibbs energy $\Gm(t,\cdot,\lambda,\theta)$ while the dissipative variable $\lambda$ is  governed by the  flow rule:
\begin{align}\label{flow-rule-physically}
\partial\zeta(\DT\lambda)=\partial_\lambda g(t,\nu,\lambda,\theta)
\end{align}
with $\partial \zeta$ denoting the subdifferential of the convex functional $\zeta(\cdot)$ and similarly $\partial_\lambda g$ is the subdifferential of the convex functional $g(t,\nu,\cdot,\theta)$. In our specific choice, \eqref{flow-rule-physically} takes the form
$\partial\delta_S^*(\DT{\lambda})+\epsilon|\DT{\lambda}|^{q-2}\DT{\lambda}+(\theta{-}\theta_{\rm c})\ccoupl\ni \varkappa\Delta^{-1}(\lambda - L \bulet \nu)$.
Furthermore, we define the specific  entropy $s$ by the standard Gibbs relation for entropy, i.e.\ $s=-g'_\theta(t,\nu,\lambda,\theta)$, and write the  entropy equation
\begin{align}\label{EntropyEquation}
\theta\DT s+{\rm div}\,j=\xi(\DT\lambda)=\,\text{ heat production rate},
\end{align}
where $j$ is the heat flux governed by the  Fourier law
\begin{align}
j=-\mathbb{K}\nabla\theta
\end{align}
with a  heat-conductivity tensor $\mathbb{K}=\mathbb{K}(\lambda,\theta)$.  In view of \eqref{form-of-zeta},
\begin{align}
\xi(\DT\lambda)=\partial\zeta(\DT\lambda){\cdot}\DT\lambda
=\delta_S^*(\DT{\lambda})+\epsilon|\DT{\lambda}|^q.
\end{align}
Now, since $s=-g'_\theta(t,\nu,\lambda,\theta)=-g'_\theta(\lambda,\theta)$,
it holds $\theta\DT s=-\theta g''_\theta(\lambda, \theta)\DT{\theta}-\theta g''_{\theta\lambda} \DT{\lambda}$. Using also $g''_{\theta\lambda}=\ccoupl$, we may reformulate the entropy equation \eqref{EntropyEquation} as the heat equation
\begin{equation}
c_\mathrm{v}(\theta) \DT{\theta}-\diver(\mathbb{K}(\lambda,\theta)\nabla \theta) =\delta_S^*(\DT{\lambda}) + \epsilon|\DT{\lambda}|^q+\ccoupl{\cdot}\theta \DT{\lambda}
\quad\text{ with }\
c_\mathrm{v}(\theta)=-\theta g ''_\theta(\theta),
\label{heatEquation0}
\end{equation}
where $c_\mathrm{v}$ is the specific heat capacity.

Altogether, we can formulate our problem  for unknowns $\theta, \nu$, and $\lambda$  which was first set and analyzed in \cite{benesova}  as
\begin{subequations} \label{ContinuousSystem}
\begin{align}
\label{ContinuousSystem-min}
&\!\!\!\!\!\!
\left.\begin{array}{llr}
\mbox{minimize} &
\displaystyle{
\int_\O\!\!\Big(\phi\bulet \nu {+} (\theta{-}\theta_{\rm c}) \ccoupl{\cdot}\lambda(t){-}\psi_0(\theta(t)) + \frac{1}{2}
m{\cdot}\nabla u_m} \\& \qquad - h(t){\cdot}m
+\frac\varkappa2\big|\nabla\Delta^{-1}(\lambda(t){-}L\bulet\nu)\big|^2
\Big)\,\d x\!\!\!\!\!\!
\\[2mm]
\mbox{subject to}\!& m=\id \bulet\nu 
\hspace{2.4em}\mbox{ on } \O,\\[2mm]
&
\displaystyle{\mbox{div} \big(\mu_0\nabla u_m-\chi_\O m \big)=0}
\hspace{3.4em}\mbox{ on } \R^d,\\[2mm] &
\nu\!\in\! {\mathscr Y}^p(\O;\R^d),\ m\!\in\! L^p(\O;\R^d),
\ u_m\!\in\! H^1(\R^d),
\end{array}\hspace*{.5em}\right\}
\text{for }t\!\in\![0,T],
\hspace*{-1em}
\\
&\partial\delta_S^*(\DT{\lambda})
+\epsilon|\DT{\lambda}|^{q-2}\DT{\lambda}+(\theta{-}\theta_{\rm c})\ccoupl\ni
\varkappa\Delta^{-1}({\rm div}\,\lambda-L\bulet \nu)
\hspace{3.2em}\text{in $Q:=[0,T]{\times}\Omega$,} \label{FlowRuleBasic}\\
&c_{\rm v}(\theta)\DT\theta-\diver(\mathbb{K}(\lambda,\theta)\nabla \theta)
= \delta_S^*(\DT{\lambda}) + \epsilon|\DT{\lambda}|^q
+\ccoupl{\cdot}\theta \DT{\lambda}
\hspace{4.5em}\text{in $Q$,} \label{heatEquation}\\
& \big(\mathbb{K}(\lambda,\theta)\nabla\theta\big){\cdot}n + b\theta=
{b}\theta_{\mathrm{ext}}\hspace{12em}
\text{ on $\Sigma:= [0,T]{\times}\Gamma$,}
\label{boundaryCondTemp}
\end{align}
\end{subequations}
where we accompanied the heat equation \eqref{EntropyEquation} by the Robin-type boundary conditions with $n$ denoting the outward unit normal to the boundary $\Gamma$, with $b\in L^\infty(\Gamma)$ a phenomenological heat-transfer coefficient, and with $\theta_{\mathrm{ext}}$ an external temperature, both assumed non-negative.
Eventually, we equip  this system with    initial conditions
\begin{align}
\lambda(0,\cdot)=\lambda_0,\qquad
\theta(0,\cdot)=\theta_0 \qquad\ \text{ on }\ \Omega,
\label{initCondi}
\end{align}

Transforming \eqref{EntropyEquation} by the  so-called enthalpy transformation, we obtain a different form of \eqref{ContinuousSystem} simpler for the analysis.  For this, let us introduce a new variable $w$, called enthalpy, by
\begin{equation}
w=\widehat c_\mathrm {v}(\theta)=\int_0^\theta c_\mathrm{v}(r) \dd r.
\label{IntroduceEnthalpy}
\end{equation}
It is natural to assume $c_\mathrm{v}$ positive, hence $\widehat c_{\rm v}$ is, for $w \geq 0$ increasing and thus invertible. Therefore, denote
 $$
\mathcalI(w):=\begin{cases}
{\widehat c_{\rm v}}^{-1}(w) & \text{if $w \geq 0$} \\
0 & \text{if $w < 0$}
\end{cases}
$$ and notice that, in the physically relevant case when $\theta \geq 0$, $\theta = \mathcalI(w)$. Thus writing the heat flux in terms of $w$ gives
\begin{align}\label{flux-in-w}
\mathbb{K}(\lambda,\theta)\nabla\theta =\mathbb{K}\big(\lambda,\mathcalI(w)\big)\nabla\mathcalI(w)
=\mathcal{K}(\lambda,w)\nabla w, \quad\text{ where }\ \mathcal{K}(\lambda,w):=
\frac{\mathbb{K}(\lambda, \mathcalI(w))}{c_\mathrm{v}(\mathcalI(w))}.
\end{align}
Moreover, the terms $(\mathcalI(w(t)){-}\theta_{\rm c}) \ccoupl{\cdot}\lambda(t)$ and $\psi_0(\theta(t))$ obviously do not play any role in the minimization \eqref{ContinuousSystem-min} and can be omitted. 
Thus we may rewrite \eqref{ContinuousSystem} in terms of $w$ as follows:
\begin{subequations} \label{ContinuousSystem1}
\begin{align}\label{ContinuousSystem1+}
&\!\!\!\!\!\!
\left.\begin{array}{llr}
\mbox{minimize} &
\displaystyle{
\int_\O\!\! \Big(\phi\bulet \nu +\frac{1}{2}
m{\cdot}\nabla u_m}- h(t){\cdot}m
+\frac\varkappa2\big|\nabla\Delta^{-1}(\lambda(t){-}L\bulet\nu)\big|^2
\Big) \,\d x\!\!\!\!\!\!
\\[2mm]
\mbox{subject to}\!& m=\id \bulet\nu,\ \
\hspace{2.4em}\mbox{ on } \O,\\[2mm]
&
\displaystyle{\mbox{div} \big(\mu_0\nabla u_m-\chi_\O m \big)=0}
\hspace{3.4em}\mbox{ on } \R^d,\\[2mm] &
\nu\!\in\! {\mathscr Y}^p(\O;\R^d),\ m\!\in\! L^p(\O;\R^d),
\ u_m\!\in\! H^1(\R^d),
\end{array}\hspace*{.5em}\right\}
\text{for }t\!\in\![0,T],
\hspace*{-1em}
\\
&\partial\delta_S^*(\DT{\lambda})
+\epsilon|\DT{\lambda}|^{q-2}\DT{\lambda}
+\big(\mathcalI(w){-}\theta_{\rm c}\big)\ccoupl\ni
\varkappa\Delta^{-1}(\lambda{-}L\bulet\nu)
\hspace{6.6em}\text{in $Q$,}
\label{FlowRuleBasic1}
\\
&\DT{w} - \diver(\mathcal{K}(\lambda, w) \nabla w)
=\delta_S^*(\DT{\lambda}) + \epsilon|\DT{\lambda}|^q
+\ccoupl \cdot \mathcalI(w)\DT{\lambda}
\hspace{6.6em}\text{in $Q$,} \label{EnthalpyEquation}
\\
&
\big(\mathcal K(\lambda,w)\nabla w\big){\cdot}n
+b\mathcalI(w)=b\theta_{\mathrm{ext}}
\hspace{13em}\text{ on $\Sigma$}.
\label{boundaryCondEnthalp}
\end{align}
\end{subequations}
Eventually, we complete this transformed system by the initial conditions
\begin{align}
\lambda(0,\cdot)=\lambda_0,\qquad
w(0,\cdot)=w_0:=\widehat c_{\rm v}(\theta_0) \qquad\ \text{ on }\ \Omega,
\label{initCond}
\end{align}
where $\lambda_0$ is the initial  phase field value, and $\theta_0$ is the initial temperature. 

Now we are ready to define a weak solution to our problem.  We denote by ${\mathscr Y}^p(\Omega;\R^d)^{[0,T]}$ the set of time-dependent Young measures, i.e., 
the set of maps $[0,T]\to{\mathscr Y}^p(\Omega;\R^d)$. We again refer to Appendix for details on Young measures.

\begin{definition}[{\sc Weak solution} \cite{benesova}]
\label{DefinOfSol}
The triple 
$$(\nu,\lambda, w){\in} ({\mathscr Y}^p(\Omega;\R^d))^{[0,T]}{\times} W^{1,q}([0,T]; L^q(\Omega;\R^{d+1})) {\times} L^1([0,T]; W^{1,1}(\Omega))$$
such that $m=\mathrm{id} \bulet \nu \in L^2(Q;\R^d)$ and $L\bulet\nu\in L^{2}(Q;\R^{d+1})$ is called a weak solution to \eqref{ContinuousSystem1}
if it satisfies:
\begin{enumerate}
\item The  minimization principle: For all $\tilde \nu$ in ${\mathscr Y}^{p}(\O;\R^d)$ and all $t \in [0,T]$
\begin{equation}
\hspace{-2em}
\Gm(t,\nu,\lambda,\mathcalI(w))\le\Gm(t,\tilde{\nu},\lambda,\mathcalI(w)).
\label{balanceLaw}
\end{equation}
\item The  magnetostatic equation:
For a.a. $t \in [0,T]$ and all $\varphi \in H^1(\R^d)$
\begin{equation}
\hspace{-2em}
\mu_0 \int_{\R^d} \nabla u_{m} \cdot \nabla \varphi\,\dd x=
\int_{\Omega} m \cdot \nabla \varphi\,\dd x.
\label{maxwellSysWeak}
\end{equation}
\item The flow rule: For any $\varphi \in L^q(Q; \R^{d+1})$
\begin{align}\nonumber
\hspace{-2em}
\int_Q\!\!\ \left( \big(\mathcalI(w){-}\theta_{\rm c}\big)\ccoupl{\cdot}
\big(\varphi{-}\DT{\lambda}\big)+\delta_S^*(\varphi)+\frac{\epsilon}{q} |\varphi|^q
+\varkappa\nabla\Delta^{-1}(\lambda{-}L\bulet\nu){\cdot}
\nabla\Delta^{-1}(\varphi{-}\DT{\lambda})
\right) \,\dd x\dd t
\\ \geq \int_{Q} \Big(\delta_S^*(\DT{\lambda})+\frac\epsilon q|\DT{\lambda}|^q \Big)\,\dd x\dd t.
\label{FlowRule}
\end{align}
\item The  enthalpy equation: For any $\varphi \in C^1(\bar{Q}), \, \varphi(T)=0$
\begin{align}
\hspace{-4em}
\int_{Q}\Big(\mathcal{K}(\lambda,w)\nabla w{\cdot}\nabla\varphi-w \DT{\varphi} \Big)\,\dd x \dd t + \int_{\Sigma} b\mathcalI(w)\varphi\,\dd S \dd t \nonumber
= \int_\Omega
w_0 \varphi(0) \,\dd x \\ + \int_{Q}\!\Big( \delta_S^*(\DT{\lambda})+\epsilon|\DT{\lambda}|^q+ \mathcalI(w)\ccoupl{\cdot}\DT{\lambda}\Big)\varphi \,\dd x \dd t + \int_{\Sigma} b\theta_{\mathrm{ext}}\varphi \,\dd S \dd t.
\label{EnthalpyEq}
\end{align}
\item The   initial conditions
in \eqref{initCond}: $\nu(0,\cdot)=\nu_0$ and $\lambda(0,\cdot)=\lambda_0$.
\end{enumerate}
\end{definition}

\noindent {\bf Data qualifications:}\\
The following  the data qualification are needed in \cite{benesova} to prove the existence of
weak solutions; cf.~\cite{benesova}:
\begin{subequations}\label{ass}
\begin{align}
\intertext{isothermal part of the anisotropy energy: $\phi \in C(\R^d)$ and}
&\exists c^A_1,c^A_2 >0,\ p>4:\ \ c^A_1(1+|\cdot|^p) \leq \phi(\cdot) \leq c^A_2(1+|\cdot|^p), \label{pGrowth}
\intertext{dissipation function: $\delta_S^* \in C(\R^{d+1})$ positively homogeneous, and}
&\exists c_{1,D},c_{2,D} >0:\ \ c_{1,D}(|\cdot|) \leq \delta_S^*(\cdot) \leq c_{2,D}(|\cdot|), \label{GrowthDissip}
\intertext{external magnetic field:}
&h \in C^1([0,T]; L^2(\Omega; \R^d)),\label{ass-f}
\intertext{specific heat capacity: $c_\mathrm{v} \in C(\R)$ and,
with $q$ from \eqref{form-of-zeta},}
&\exists c_{1,\theta},c_{2,\theta}>0,\ \omega_1 \geq \omega \geq q',\
c_{1,\theta} (1{+}\theta)^{\omega-1} \leq c_\mathrm{v}(\theta) \leq c_{2,\theta}(1{+}\theta)^{\omega_1-1},
\label{OmegaGrowth}
\intertext{heat conduction tensor: $\mathcal{K}\in C(\R^{d+1}\times\R;\R^{d \times d})$ and }
&\exists C_K,\kappa_0 >0\ \forall\chi\in\R^d: \mathcal{K}(\cdot , \cdot) \leq C_K, \, \, \, \, \chi^\mathrm{T}\mathcal{K}(\cdot, \cdot) \chi \geq \kappa_0 |\chi|^2,
\label{Sub6Growth}
\intertext{external temperature:}
&\theta_\mathrm{ext} \in L^1(\Sigma),\ \ \theta_\mathrm{ext}\ge0,
\ \text{ and }\ b \in L^\infty(\Sigma),\ \ b\ge0,
\label{ExternalTemp}
\intertext{initial conditions:}
&\nu_0\in{\mathscr Y}^p(\Omega;\R^d) \text{ solving \eqref{balanceLaw} },\quad \lambda_0\in L^q(\O;\R^{d+1}),
\quad w_0=\widehat c_\mathrm{v}(\theta_0)\in L^1(\O)
\text{ with }\theta_0\ge0. \label{ass-IC}
\end{align}
\end{subequations}

The following theorem is proved in \cite{benesova}.

\begin{theorem}\label{thm-main}
Let \eqref{ass} hold. Then at least one weak solution $(\nu,\lambda, w)$
to the problem \eqref{ContinuousSystem1}
in accord with Definition~\ref{DefinOfSol} does exist. Moreover, some of
these solutions satisfies also
\begin{align}\label{estimates-of-w}
&w\in L^r([0,T]; W^{1,r}(\Omega))\,\cap\,W^{1,1}(I;W^{1,\infty}(\Omega)^*)
\qquad\text{ with }1\le r<\frac{d{+}2}{d{+}1}.
\end{align}
\end{theorem}

The proof of the Theorem~\ref{thm-main}  in \cite{benesova} exploits the following  time-discrete approximations  which  also create  basis for our fully discrete solution. Given $ T>0$ and $T/\tau\in\N$
we call the triple $(\nu_\tau^k, \lambda_\tau ^k, w_\tau ^k) \in {\mathscr Y}^{p}(\O;\R^d) \times L^{2q}(\Omega; \R^{d+1})\times H^1(\Omega)$ the \emph{discrete weak solution} of \eqref{ContinuousSystem1} subject to boundary condition \eqref{boundaryCondEnthalp} at time-level $k$, $k = 1 \ldots, T/\tau$, if it satisfies:
\begin{subequations}\label{def-disc}
\begin{enumerate}
\item The time-incremental {\bf minimization problem} with given $\lambda_\tau ^{k-1}$ and $w_\tau ^{k-1}$:
\begin{align}
\left.
\hspace{-2ex}
\begin{array}{ll}
\mathrm{Minimize}\!\!& \displaystyle{\Gm(k\tau,\nu,\lambda,\mathcalI(w_\tau^{k-1}))
+\tau\!\int_\Omega\!\left(|\lambda|^{2q}
+\delta_S^*\Big(\frac{\lambda{-}\lambda_\tau ^{k-1}}{\tau}\Big)
+\frac{\epsilon}{q}\Big|\frac{\lambda{-}\lambda_\tau ^{k-1}}{\tau}\Big|^q \right) \dd x}
\\[.7em]
\mathrm{subject \, to}\!\!&
(\nu,\lambda)\in {\mathscr Y}^{p}(\O;\R^d) \times L^{2q}(\Omega; \R^{d+1}).
\end{array}\!\right\}\!\!\!
\label{BalanceEqDis}
\end{align}
with $\Gm$ from \eqref{def-of-psi}.
\item \reviewMK{The {\bf Poisson problem}:}
For all $\varphi \in H^1(\R^d)$
\begin{equation}
\int_{\R^d} \nabla u_{m_\tau^k}{\cdot}\nabla \varphi\,\dd x = \int_{\Omega}m_\tau^k {\cdot}\nabla \varphi\,\dd x\qquad\text{ with }\ m_\tau^k=\mathrm{id} \bulet \nu_\tau^k.
\label{MaxwellDis}
\end{equation}
\item The {\bf enthalpy equation}: For all $\varphi \in H^1(\Omega)$
\begin{align}
\int_{\Omega}\left(\frac{w_\tau^k{-}w_\tau ^{k-1}}{\tau}\varphi + \mathcal{K}(\lambda_\tau^k, w_\tau^k) \nabla w_\tau ^k{\cdot}\nabla\varphi \right)\,\dd x
+ \int_{\Gamma} b^k_\tau\mathcalI(w_\tau ^k)\varphi\,\dd S \nonumber
= \int_{\Gamma} b^k_\tau\theta^k_{\mathrm{ext}, \tau}\varphi\,\dd S
\\
+\int_{\Omega}\left(\delta_S^*
\Big(\frac{\lambda_\tau^k{-}\lambda_\tau^{k-1}}{\tau}\Big)
+\epsilon\Big|\frac{\lambda_\tau^k{-}\lambda_\tau^{k-1}}{\tau}\Big|^q
\mathcalI(w_\tau^k)\ccoupl{\cdot}
\frac{\lambda^k{-}\lambda^{k-1}}{\tau} \right)\varphi\,\dd x.
\label{EnthalpyEqDis}
\end{align}
\item For $k=0$ the {\bf initial conditions} in the following sense
\begin{align}\label{initCondDis}
\nu_\tau ^0=\nu_0,\qquad \lambda_\tau ^0=\lambda_{0,\tau},
\qquad w_\tau ^0=w_{0,\tau}\ \ \ \text{ on }\O.
\end{align}
\end{enumerate}
\end{subequations}

In \eqref{initCondDis}, we denoted by $\lambda_{0,\tau}\in L^{2q}(\O;\R^{d+1})$
and $w_{0,\tau}\in L^2(\O)$
respectively suitable approximation of the
original initial conditions $\lambda_0\in L^q(\O;\R^{d+1})$ and
$w_0\in L^1(\O)$ such that
\begin{subequations}
\begin{align}\label{lambdaZeroApprox}
&\lambda_{0,\tau}\to\lambda_0\ \text{ strongly in $L^q(\O;\R^{d+1})$, and }\
\|\lambda_{0,\tau}\|_{L^{2q}(\O;\R^{d+1})} \leq C \tau^{-1/(2q+1)},
\\
&w_{0,\tau}\to w_0\ \text{ strongly in $L^1(\O)$, and }\ w_{0,\tau}\in L^2(\O).
\end{align}
\end{subequations}
Moreover $\theta^k_{\mathrm{ext}, \tau} \in L^2(\Gamma)$ and $b^k_\tau \in L^\infty (\Gamma)$ are defined in such a way that their piecewise constant interpolants
\begin{equation*}
\big[\bar\theta_{\mathrm{ext},\tau},\bar b_\tau](t):=\big(\theta^k_{\mathrm{ext}, \tau},b^k_\tau,)
\qquad\qquad\text{ for $\ (k{-}1)\tau<t\le k\tau$,\ \ $k=1,...,K_\tau$}.
\end{equation*}
satisfy
\begin{equation}
\bar\theta_{\mathrm{ext}, \tau} \to \theta_\mathrm{ext}\ \text{ strongly in $L^1(\Sigma)$ and }\ \bar b_\tau \stackrel{*}{\rightharpoonup} b \ \text{ weakly* in $L^\infty(\Sigma)$.}
\label{ConvTempData}
\end{equation}

We introduce the notion of \emph{piecewise affine} interpolants
$\lambda_\tau$ and $w_\tau$ defined by
\begin{align}\nonumber
\big[\lambda_\tau,w_\tau\big](t):=\frac{t-(k{-}1)\tau}\tau \big(\lambda_\tau^k,w_\tau^k\big) +\frac{k\tau-t}\tau \big(\lambda_\tau^{k-1},w_\tau^{k-1}\big)
\qquad\text{ for $\ t\in[(k{-}1)\tau,k\tau]\ $}
\end{align}
with $\ k=1,...,T/\tau$.
In addition, we define the backward \emph{piecewise constant interpolants}
$\bar{\nu}_\tau$, $\bar{\lambda}_\tau$, and $\bar{w}_\tau$ by
\begin{align}
\big[\bar{\nu}_\tau,\bar{\lambda}_\tau,\bar{w}_\tau\big](t):=\big(\nu_\tau^k,\lambda_\tau^k,w_\tau^k\big)
\qquad\qquad\text{ for $\ (k{-}1)\tau<t\le k\tau$,\ \ $k=1,...,T/\tau$}.
\end{align}
Finally, we  also need the piecewise constant interpolants of delayed  enthalpy and magnetization  $\underline{w}_\tau$, $\underline{m}_\tau$ defined by
\begin{align}\label{w-backward}
[\underline{w}_\tau(t), \underline{m}_\tau(t)]:= [w_\tau^{k-1}, \id \bulet \nu_\tau^{k-1}] \qquad\qquad \text{ for $\ (k{-}1)\tau<t\le k\tau$,\ \ $k=1,...,T/\tau$}.
\end{align}

\subsection{Energetics} 
In this section we summarize  some basic energetic estimates available for our model.
First we define   the purely magnetic part of the Gibbs free energy $\mathfrak{G}$  as
\begin{equation}\label{def-of-G-mag}
\mathfrak{G}(t,\nu, \lambda):=\int_\Omega \phi \bulet \nu -
h(t){\cdot}m\,\dd x + \int_{\R^d} \frac{1}{2}|\nabla u_{m}|^2\,\dd x
+\frac\varkappa2\big\|\lambda {-} L \bulet \nu\big\|_{H^{-1}(\Omega; \R^{d+1})}^2.
\end{equation}

The purely magnetic part of the Gibbs energy satisfies (see \cite[Formula (4.19)]{benesova}) the following energy inequality 
\begin{equation}\label{magneticbalance}
\mathfrak{G}(t_\ell, \bar{\nu}_\tau (t_\ell), \bar{\lambda}_\tau(t_\ell)) \leq  \mathfrak{G}(0, \bar{\nu}_\tau (0), \bar{\lambda}_\tau(0)) + \int_0^{t_\ell} \Big( \int_\Omega \DT{h}_\tau \cdot \bar{m}_\tau \dd x + \varkappa \llangle \bar{\lambda}_\tau - L \bulet \bar{\nu}_\tau, \DT{\lambda}_\tau\rrangle \Big) \dd t
\end{equation}
with $t_\ell = \ell \tau$. 

As $(\nu_\tau^k,\lambda_\tau^k)$ is a minimizer of \eqref{BalanceEqDis}, the partial sub-differential of the cost functional with respect to $\lambda$ has to be zero at $\lambda_\tau^k$. This condition holds at each time level and, thus,  summing up for $k=0,\ldots,\ell$ gives 
\begin{align}
&\!\!\int_0^{t_\ell}\!\!\int_\Omega \Big(\delta_S^*( \DT{
\lambda}_\tau)+\frac{\epsilon}{q}|\DT{\lambda}_\tau|^q \Big) \dd x \dd t
\leq \int_0^{t_\ell}\!\!\bigg(\varkappa\llangle \bar{\lambda}_\tau{-}L \bulet\bar{\nu}_\tau , v_\tau{-}\DT{\lambda}_\tau \rrangle
\nonumber \\&\quad+ \int_\Omega \Big(\big(\mathcalI(\underline{w}_\tau)-\theta_{\rm c}\big) \ccoupl{\cdot}(v_\tau{-}\DT{\lambda}_\tau) + 2q\tau |\bar{\lambda}_\tau|^{2q-2}\bar{\lambda}_\tau(v_\tau{-}\DT{\lambda}_\tau)
+\delta_S^*(v_\tau)+\frac{\epsilon}{q} |v_\tau|^q \Big)\dd x\bigg) \dd t
\label{FlowRuleDisI},
\end{align}
where $v_\tau$ is an arbitrary test function such that $v_\tau(\cdot, x)$ is piecewise constant on the intervals $(t_{j-1}, t_j]$ and $v_\tau(t_j,\cdot) \in L^{2q}(\Omega; \R^{d+1})$ for every $j$.

Hence, for $v_\tau=0$ we get the energy balance of the thermal part of the Gibbs energy, namely 
\begin{align}\label{thermalbalance}
&\!\!\int_0^{t_\ell}\!\!\int_\Omega \Big(\delta_S^*( \DT{
\lambda}_\tau)+\frac{\epsilon}{q}|\DT{\lambda}_\tau|^q \Big) \dd x \dd t
\leq \int_0^{t_\ell}\!\!\bigg(-\varkappa\llangle \bar{\lambda}_\tau{-}L \bulet\bar{\nu}_\tau , \DT{\lambda}_\tau \rrangle
- \int_\Omega \big(\mathcalI(\underline{w}_\tau)-\theta_{\rm c}\big) \ccoupl{\cdot}\DT{\lambda}_\tau + 
2q\tau |\bar{\lambda}_\tau|^{2q-2}\bar{\lambda}_\tau\DT{\lambda}_\tau \bigg)\dd t\ .
\end{align}

This inequality couples the dissipated energy and temperature evolution.

\section{Numerical approximations and computational examples}\label{numerics}

Dealing with a numerical solution, we have to find  suitable spatial approximations for $\nu$, $u_m$, $w$, and $\lambda$ in each time step. In our numerical method, we require that \eqref{moment-constraint} is  satisfied which means that knowing the Young measure $\nu$ we can easily calculate the momenta $\lambda$. We present a spatial discretization of involved quantities in each time step.  

 The domain $\O$ of the ferromagnetic body is discretized by a regular triangulation $\mathcal T_\ell$ in triangles (in 2D) or in tetrahedra (in 3D) for $\ell\in\N$ which will be called elements. The triangulations are nested, i.e.,  that $\mathcal{T}_\ell\subset\mathcal{T}_{\ell+1}$, so that the discretizations are finer as $\ell$ increases. Let us now describe the approximation.

{\it Young measure.}  Young measures are parametrized (by $x\in\O)$  probability measures supported on $\R^d$.  Hence, we need to handle their discretization in $\O$ as well as in $\R^d$.  Our aim is to approximate a general Young measure by  a  convex combination of a finite number of  Dirac measures (atoms) supported on $\R^d$ such that this convex combination is elementwise constant.  Let us now describe a rigorous procedure how to achieve this goal. We first omit the time discretization parameter $\tau$ and discuss the  discretization of the Young measure in $\O$.  In order to approximate a Young measure $\nu$,  we follow \cite{carstensen-roubicek, matache} and define for  $z\in L^\infty(\O)\otimes C^p(\R^d)$  the following projection operator ($\mathcal{L}^d$ denotes the $d$-dimensional Lebesgue measure) 
$$
[\Pi^1_\ell z](x,s)=\frac{1}{\mathcal{L}^d(\triangle)}\int_\triangle z(\tilde x,s)\,\d \tilde x\ \,  \mbox{ if } x\in\triangle\in\mathcal{T}_\ell\ .
$$
Notice that $\Pi^1_\ell$ is elementwise constant in the $x$-variable.  We now turn to a discretization of $\R^d$ in terms of   large cubes in $\R^d$, i.e.,   for $\alpha\in\N$  we consider a cube $B_\alpha:=[-\alpha,\alpha]^d$ (i.e. we call it ``a cube'' even if $d=2$) which is discretized into  $(2\alpha/n)^d$ smaller cubes with the edge length $2\alpha/n$ for some $n\in\N$. Corners of small cubes are called nodal points.  We  define $Q_1$ elements on the  cube $B_\alpha\in\R^d$ which consist of  tensorial products of affine functions in each spatial variable of $\R^d$. In this way, we find  basis functions $f_i:B_\alpha\to\R$ for $i=1,\ldots,(n+1)^d$ such that $f_i\ge 0$ and $\sum_{i=1}^{(n+1)^d} f_i(s)=1$ for all $s\in\R^d$. Moreover, if $s_j$ is the $j$-th nodal point then  $f_i(s_j)=\delta_{ij}$, where $\delta_{ij}$ is the Kronecker symbol. 
Further, each $f_i$ can be continuously  extended to $\R^d\setminus B_\alpha$  and such an  extended  function can even vanish at infinity, i.e., it belongs to $C_0(\R^d)$.
This construction defines a projector $L^\infty(\O)\otimes C^p(\R^d)\to L^\infty(\O)\otimes C^p(\R^d) $ as 
$$
[\Pi^2_{\alpha,n}z](x,s):=\sum_{i=1}^{(n+1)^d} z(x,s_i)f_i(s)\ .
$$
    
Finally, we define $\Pi_{\ell,\alpha,n}:=\Pi^1_\ell\circ \Pi^2_{\alpha,n}$, so that 

$$
[\Pi_{\ell,\alpha,n}z](x,s):= \frac{1}{\mathcal{L}^n(\triangle)}\sum_{i=1}^{(n+1)^d}\int_\triangle z(\tilde x,s_i)v_i(s)\,\d \tilde x\ \,  \mbox{ if } x\in\triangle\in\mathcal{T}_\ell\ .
$$
If we now  take $\nu\in\mathscr{Y}^p(\O;\R^d)$ and denote $l:=(\ell,\alpha,n)$   we calculate 
\begin{equation}\label{projectedmeasure}
\int_\O\int_{\R^d}[\Pi_{l}z](x,s)\nu_x(\d s)\,\d x=\int_\O\int_{\R^{d}} z(x,s)[\nu_l]_x(\d s )\,\d x\ ,
\end{equation}
where for $x\in\O$
\begin{equation}\label{discreteYm}
[\nu_l]_x:=\sum_{i=1}^{(n+1)^d} \xi_{i,l}(x)\delta_{s_i}\ ,
\end{equation}
with
$$\xi_{i,l}(x):= \frac{1}{\mathcal{L}^d(\triangle)}\int_\triangle \int_{\R^d}f_i(s)\nu_x(d s)\,\d x\ ,\ x\in\triangle\in\mathcal{T}_\ell\ .$$ Let us denote 
the subset of Young measures from $\mathscr{Y}^p(\O;\R^d)$ which are in the form of \eqref{discreteYm} by $\mathscr{Y}_l^p(\O;\R^d)$. 
Notice that $\xi_{i,l}\ge 0$ and that $\sum_{i=1}^{(n+1)^d}\xi_{i,l}=1$. Hence, the projector $\Pi_l$ corresponds to approximation of $\nu$ by a spatially piecewise constant Young measure which can be written as a convex combination of Dirac measures (atoms). We refer to \cite{roubicek0} for a thorough description of various kinds of Young measure approximations.   In order to indicate that the  measure is time-dependent we write in the $k$-th time-step 
$$
[\nu_{l, \tau}^k]_x:=\sum_{i=1}^{(n+1)^d} \xi^k_{i,l,\tau}(x)\delta_{s_i}\ .
$$ 

\begin{figure}
\center
\begin{minipage}{0.45\textwidth}
\center
\includegraphics[width=\textwidth]{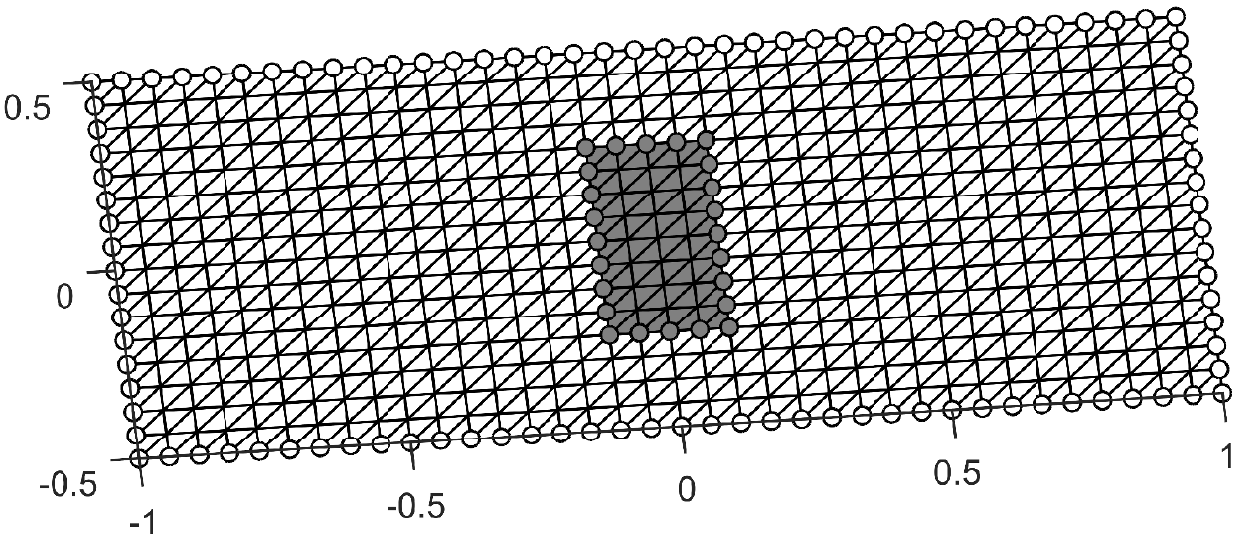}
\end{minipage}
\begin{minipage}{0.1\textwidth}
\end{minipage}
\begin{minipage}{0.45\textwidth}
\center
\includegraphics[width=\textwidth]{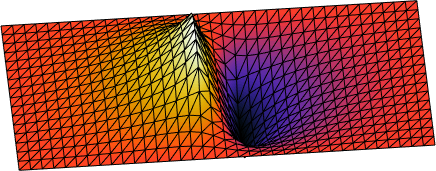}
\end{minipage}
\caption{Example of the outer triangulation $\hat{ \mathcal T}$ containing the magnet body triangulation $\mathcal T$ (in gray) is shown in the left. The right part displays an example of the magnetostatic potential $u_{m}$ approximated as the scalar nodal and elementwise linear function (P1 elements function) satisfying zero Dirichlet condition in the boundary nodes of  $\hat{ \mathcal T}$.}
\label{figure:magnet_and_magnetostatic_potential}
\end{figure}

{\it Magnetostatic potential.} Following \cite{carstensen-prohl}, we simplify the calculation of the  reduced Maxwell system  in   magnetostatics by assuming  that the magnetostatic potential $u$ vanishes outside  a large bounded domain $\hat\Omega\supset  \Omega$. Hence, given $m\in L^p(\O;\R^d)$, we solve the Poisson problem ${\rm div} (\mu_0\nabla u_m)={\rm div} (\chi_\O m) $ on $\hat\O$ with homogeneous Dirichlet boundary condition $u_m=0$ on $\partial\hat\O$.  The set $\hat\Omega$ is discretized by an outer triangulation $\hat{ \mathcal T}_\ell$ that contains the triangulation $\mathcal{T}_\ell$ of the ferromagnetic magnetic body. Then, the magnetostatic potential 
\begin{equation}
u_{m_{l, \tau}^k} \in P^1_0(\hat{ \mathcal T_\ell})
\label{space_magnetostatic_potential}
\end{equation}
in the $k$-th time-step  is approximated in the space $P^1_0(\hat{ \mathcal T_\ell})$ of scalar nodal and elementwise linear functions defined on the triangulation $ \hat{ \mathcal T_\ell}$ and  satisfying zero Dirichlet boundary conditions on the triangulation boundary $\partial \hat{ \mathcal T_\ell}$ . For illustration, see Figure \ref{figure:magnet_and_magnetostatic_potential}.
The magnetization vector 
\begin{equation}
m_{l,\tau}^k \in P^0(\mathcal T_\ell)^d
\label{space_magnetization}
\end{equation}
in the $k$-th time-step  is approximated in the space $P^0(\mathcal{T}_\ell)^d$ of vector and elementwise constant functions.  Another numerical approaches to solutions of magnetostatics using e.g.~BEM are also available \cite{andjelic-of-steinbach-urtaler}. 

{\it Enthalpy.} The enthalpy 
\begin{equation}
w_{\ell,\tau}^k \in P^1(\mathcal T_\ell)
\label{space_enthalpy}
\end{equation}
in the $k$-th time-step  is approximated in the space $P^1(\mathcal T_\ell)$ of scalar nodal and elementwise linear functions.

Having time and spatial discretizations we can set up an algorithm to solve the problem which is just \eqref{def-disc} with additional spatial discretization.   Finally, we apply the spatial discretization just described and we arrive at the following problem. \\

Given  spatially discretized  boundary condition \eqref{boundaryCondEnthalp} and   $k = 1, \ldots , T/\tau$ we  solve:
\begin{subequations}\label{def-disc1}
\begin{enumerate}
\item The  minimization problem with given $w^{k-1}_{\ell,\tau}\in P^1(\mathcal T_\ell)^d $  with $\lambda^{k-1}_{l,\tau}:=L\bulet\nu^{k-1}_{l,\tau}$:
\begin{align}
\left.
\begin{array}{ll}
\mathrm{Minimize}\!\!& \displaystyle{\Gm(k\tau,\nu,\lambda,\mathcalI(w_{\ell,\tau}^{k-1}))
+\tau\!\int_\Omega\!\left(|\lambda|^{2q}
+\delta_S^*\Big(\frac{\lambda{-}\lambda_{l,\tau}^{k-1}}{\tau}\Big)
+\frac{\epsilon}{q}\Big|\frac{\lambda{-}\lambda_{l,\tau}^{k-1}}{\tau}\Big|^q \right) \dd x}
\\[.7em]
\mathrm{subject \, to}\!\!&
\nu\in {\mathscr Y}_l^{p}(\O;\R^d)\ ,\ \lambda:=L\bulet\nu
\end{array}\!\right\}\!\!\!
\label{BalanceEqDis1}
\end{align}
with $\Gm$ from \eqref{def-of-psi}.
\item The Poisson problem:
For all $v \in P_0^1(\hat{\mathcal T}_\ell)$
\begin{equation}
\mu_0 \int_{\R^d} \nabla u_{m_{l,\tau}^k}{\cdot}\nabla \varphi\,\dd x = \int_{\Omega}m_{l,\tau}^k {\cdot}\nabla \varphi\,\dd x\qquad\text{ with }\ m_{l,\tau}^k=\mathrm{id} \bulet \nu_{l,\tau}^k.
\label{MaxwellDis1}
\end{equation}
\item The enthalpy equation: For all $\varphi \in P^1(\mathcal T_\ell)$
\begin{align}
\int_{\Omega}\left(\frac{w^k_{\ell,\tau}{-}w_{\ell,\tau}^{k-1}}{\tau}\varphi + \mathcal{K}(\lambda^k_{l,\tau}, w^k_{\ell,\tau}) \nabla w_{\ell,\tau} ^k{\cdot}\nabla\varphi \right)\,\dd x
+ \int_{\Gamma} b\mathcalI(w_{\ell,\tau} ^k)\varphi\,\dd S \nonumber
= \int_{\Gamma} b\theta^k_{\mathrm{ext}, \tau}\varphi\,\dd S
\\
+\int_{\Omega}\left(\delta_S^*
\Big(\frac{\lambda^k_{l,\tau}{-}\lambda_{l,\tau}^{k-1}}{\tau}\Big)
+\epsilon\Big|\frac{\lambda^k_{l,\tau}{-}\lambda_{l,\tau}^{k-1}}{\tau}\Big|^q
+\mathcalI(w^k_{\ell,\tau})\ccoupl{\cdot}
\frac{\lambda^k_{l,\tau}{-}\lambda^{k-1}_{l,\tau}}{\tau} \right)\varphi\,\dd x.
\label{EnthalpyEqDis1}
\end{align}
\item For $k=0$ the  initial conditions:
\begin{align}\label{initCondDis0}
\lambda_{l,\tau} ^0=\lambda_{0,l},\qquad \qquad w_{\ell,\tau} ^0=w_{0,\ell}\ \ \ \text{ on }\O \ ,
\end{align}
\end{enumerate}
\end{subequations}
where $\lambda_{0,\ell}=L\bulet\nu_{0,\ell}$ is  calculated via \eqref{projectedmeasure} and $w_{0,\ell}$ is a piecewise affine approximation of $w_0$. There is no initial condition for $\lambda_{\ell,\tau} ^0$ as it is now fully determined by $\nu_{0,\ell}$. \\
 \begin{figure}
\center
\begin{minipage}{0.45\textwidth}
\center
\includegraphics[height=0.3\textheight]{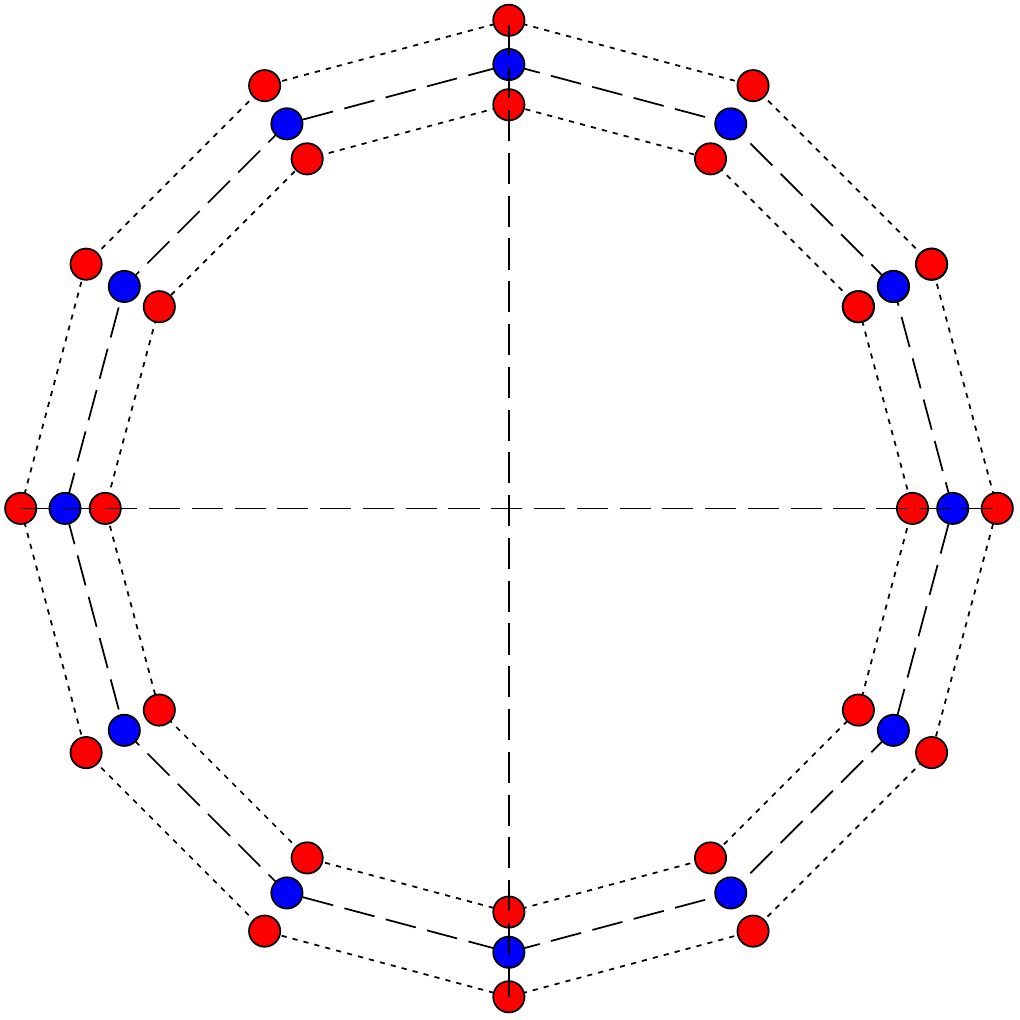}
\end{minipage}
\begin{minipage}{0.1\textwidth}
\end{minipage}
\begin{minipage}{0.45\textwidth}
\center
\includegraphics[height=0.3\textheight]{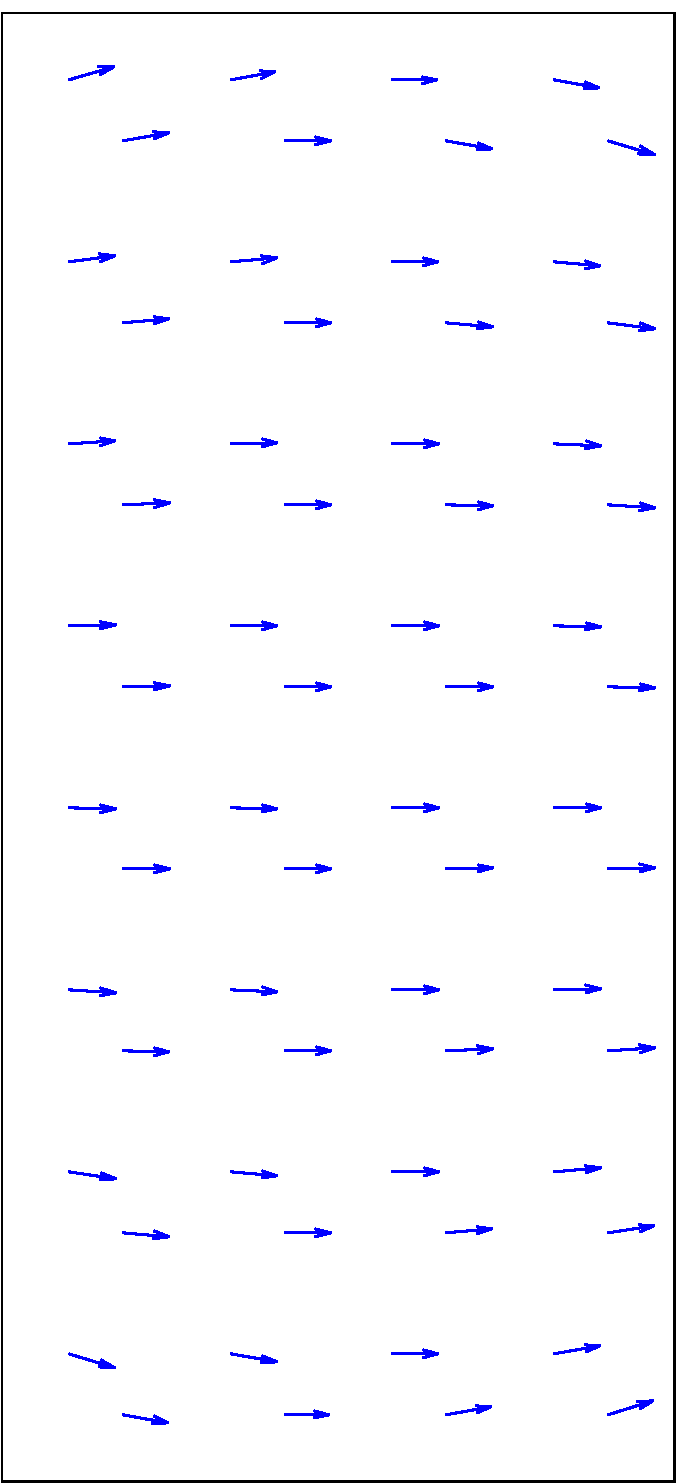}
\end{minipage}
\caption{An example of uniformly distributed Dirac atoms on the left: Each atom is specified by its angle $\varphi_i$ and radius $r_i $  for $i=1, \dots, N$. Here, $N=36$ and Dirac atoms are placed on ``the main sphere'' with radius $1$ (blue colored atoms in the color scale or dark colored atoms in the gray scale) and additional two spheres with radii $\frac{1}{1.1}$ and $1.1$ (red colored atoms in the color scale or gray colored atoms in the gray scale). An example of magnetization $m$ is displayed on the right. Each vector (arrow) corresponds to value of $m$ in one element and its orientation is given as a convex combination of Dirac atoms multiplied by the value of $ p_\tau^k$, see \eqref{convex_combination}.}
\label{figure:meassure_supports_and_magnetization}
\end{figure}

In computations, several simplifications were taken to account. First of all,  we assume 
\begin{equation} 
d=2, \qquad q=2.
\end{equation} 
In view of \eqref{moment-constraint}, the macroscopic magnetization $m$ is elementwise constant and it is the first moment of $\nu_l$. As the anisotropy energy density is minimized for a given temperature  on a sphere in $\R^d$ we put the support of the Young measure $\nu_l$ on this  sphere and its vicinity to decrease the number of variables in our problem. In what follows, the number of Dirac atoms in $\nu_l$ is denoted by $N\in\N$. It is then  convenient to work in polar coordinates  where $r_i$ is the radius and $\varphi_i$ the corresponding angle of the $i$-th atom. Hence, we have   
\begin{equation}
m_{l,\tau}^k=\lambda_{1,l,\tau}^k=p_\tau^k \sum_{i=1}^N \xi_{i,l,\tau}^k r_i \, (\cos(\varphi_i), \sin(\varphi_i)), 
\qquad \qquad
\lambda_{2,l,\tau}^k= (p_\tau^k)^2 \sum_{i=1}^N \xi_{i,l,\tau}^k r_i^2, 
\qquad \qquad \sum_{i=1}^N \xi_{i,l,\tau}^k  = 1, 
\label{convex_combination}
\end{equation}
where coefficients $\xi_{i,l,\tau}^k \in [0,1], i=1, \dots, N$, and $p_\tau^k$ depends on temperature in the following way:
$$
p^k_\tau(\theta):=\begin{cases}
\sqrt{(\theta_{\rm c}-\theta)a_0/(2b_0)} & \text{ if $\theta_{\rm c}>\theta$,}\\
p_{\rm{par}} &\text{ otherwise.}
\end{cases}
$$
A small parameter $p_{\rm{par}}>0$ is introduced   which allows for nonzero magnetization and   increase of the temperature due to the change of magnetization even in the paramagnetic mode.  The number $N$ and values of radii $r_i$ and angles $\varphi_i$ are given a priori and influence possible directions of magnetization, see Figure \ref{figure:meassure_supports_and_magnetization}. The coefficients of the convex combinations and $p_\tau^k$  in the $k$-th time-step}
\begin{equation}
\xi_{i,l,\tau}^k, p_\tau^k \in P^0(\mathcal T_\ell) 
\label{space_meassure_support}
\end{equation}
for all $i=1, \dots, N$ are approximated in the space $P^0(\mathcal T_\ell)$ of scalar and elementwise constant functions. 
We assume that for $H_{\rm c},h_{\rm c}>0$
 $$
S:=\{\lambda=(\lambda_1,\lambda_2)\in\R^2\times\R:\, |\lambda_1|\le H_{\rm c}\,\&\, |\lambda_2|\le h_{\rm c}\}\ .$$ Then for $\eta\in\R^2\times\R$
\begin{equation}\label{norm-term}\delta_S^*(\eta)=\max_{\lambda\in S} \eta\cdot\lambda=H_{\rm c} |\eta_1|+h_{\rm c}|\eta_2|\ . 
\end{equation} 
where $H_c$ represents the coercive force of the magnetic material. Then the minimization problem \eqref{BalanceEqDis1} can be expressed in unknown coefficients $\xi_{i,l,\tau}^k$, $i=1, \dots, N$ only. The functional in \eqref{BalanceEqDis1} contains a nondifferentiable norm term \eqref{norm-term}, and its evaluation requires to solve the magnetostatic potential $u_{m_{l,\tau}^k}$ from  the Poisson problem \eqref{MaxwellDis1} with zero boundary conditions. The size of the matrix in  the discretized Poisson problem equals to the number of free nodes in the triangulation $\hat{ \mathcal T}_\ell$. After coefficients $\xi_{i,l,\tau}^k$ for $i=1, \dots, N$ are computed, the enthalpy $w^k_{\ell,\tau}$ is solved from the enthalpy equation \eqref{EnthalpyEqDis1}. We consider the case 
\begin{equation}
\mathbb{K}(\lambda,\theta)={\rm const.}, \qquad c_v(\theta)={\rm const.}
\end{equation}
of the constant heat-conductivity $\mathbb{K}$ and the constant heat capacity $c_v$.  Therefore, the enthalpy equation \eqref{EnthalpyEqDis1} can be discretized as a linear system of equations combining stiffness and mass matrices from the discretization of a second order elliptic partial differential equation using $P^1$ elements. Therefore, the size of both matrices is equal to the number of all nodes in the triangulation $\mathcal T_\ell$.

As an example of computation, we consider a large domain $\hat \Omega$ and a magnet domain $\Omega$, where
$$\hat \Omega = (-1,1 ) \times (-\frac{1}{2}, \frac{1}{2}), \qquad \Omega=(-\frac{1}{9}, \frac{1}{9}) \times (-\frac{1}{4}, \frac{1}{4})$$
with a triangulation shown in Figure \ref{figure:magnet_and_magnetostatic_potential} (left). A Young measure was discretized using 36  Dirac measures grouped in three spherical layers as shown in Figure \ref{figure:meassure_supports_and_magnetization} (left). 

Physical parameters   were chosen  to show qualitative results only and they obviously do not correspond to any realistic material. We consider
\begin{itemize}
\item $\phi_\mathrm{poles}(m)=m_1^2$, where $m=(m_1,m_2)$ and  $m$ is measured in A$/$m,
\item the coercive force $H_c= 100 \,{\rm T} $ - this value provides a hysteresis width visible in all figures,
\item $h_{\rm c}=1 {\rm Tm}/{\rm A}$
\item $p_{\rm{par}} =0.1 $
\item the parameter\footnote{$\epsilon$ stands in front of $\lambda$ whose units depend on a particular component. Hence, to avoid constants of value one which only carry  SI units we do not specify the unit of $\epsilon$.}  $\epsilon=10^{-6}$
\item the initial temperature inside magnet $\theta_0=1300 \,{\rm K}$,  the Curie temperature $\theta_c= 1388$ K and the constant external temperature around the magnet body is $\theta_\mathrm{ext}=1100\,{\rm K} $, 
\item the coefficient  $b=0.001 {\rm  W}/({\rm mK}$)  in the Robin-type boundary condition, the heat conductivity coefficient ($\mathbb{I}$ stands for the identity matrix in $\R^{2\times 2}$) $\mathbb{K} =100\, \mathbb{I} {\rm W}/{\rm m\,K}$ and the heat capacity $c_v=420 {\rm J}/({\rm m}^3{\rm K})$,
\item the coefficients in the thermo-magnetic coupling $a_0=1\, {\rm J}/({\rm KmA}^2)$, $b_0=1\, {\rm Jm}/{\rm A}^4$,
\item the uniaxial cyclic magnetic field $h(t)=3 H_c(h_x(t),0)${\rm  T}, where $t=0,\dots,80$ and $h_x$ is a cyclic periodic function with the period $10$ and the amplitude $1$. 
\end{itemize}

\begin{figure}
\center
\includegraphics[width=\textwidth,height=0.25\textheight]{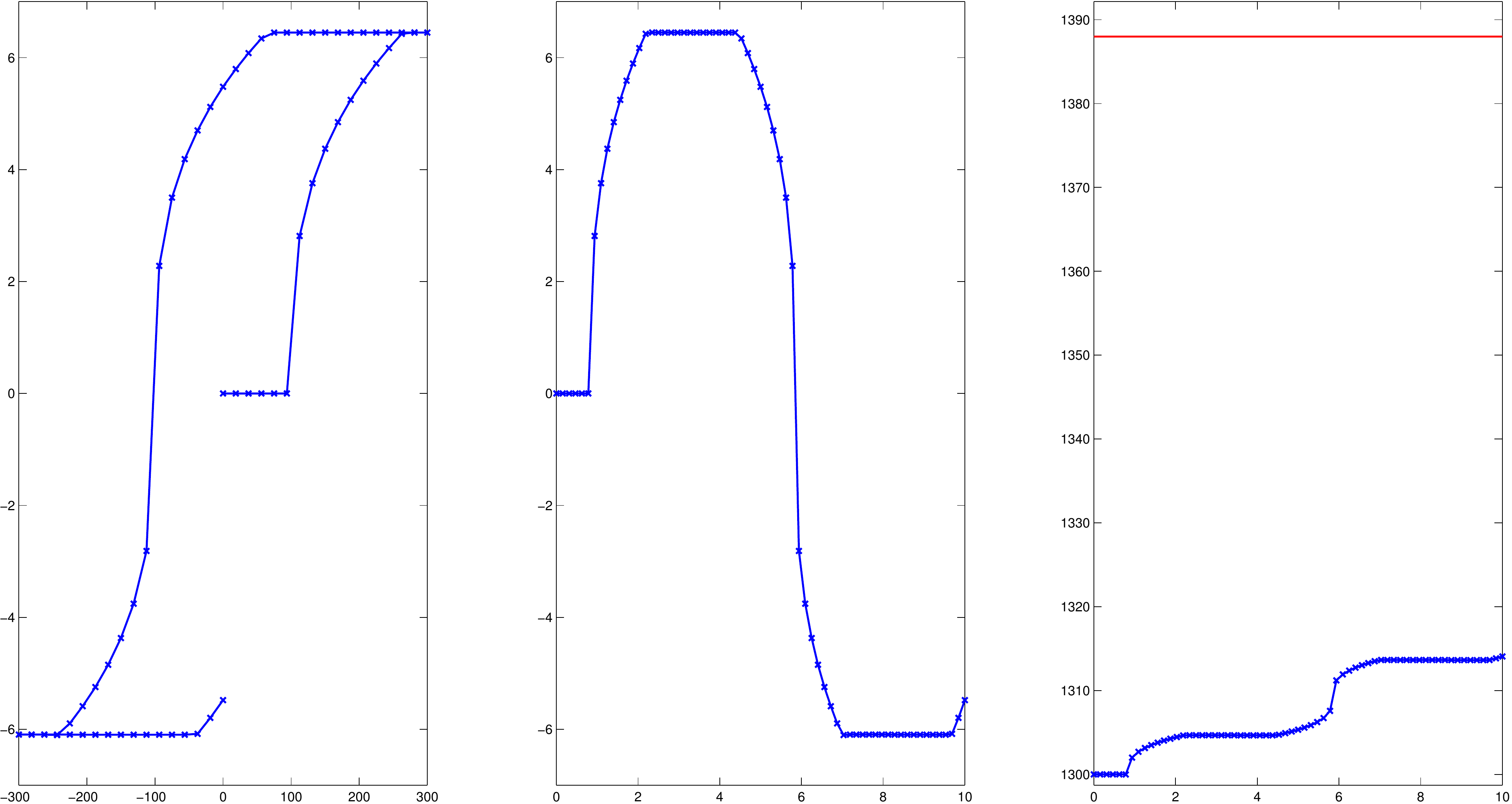}
\vspace{-0.7cm}
\caption{Average values of fields after one cycle of external forces: magnetization in x-direction versus external field (left), 
magnetization in x-direction versus time (middle), temperature versus time (right) never reaching  the Curie temperature indicated by the red horizontal line.} 
\label{figure-one_loop}
\vspace{0.5cm}
\includegraphics[width=\textwidth,height=0.25\textheight]{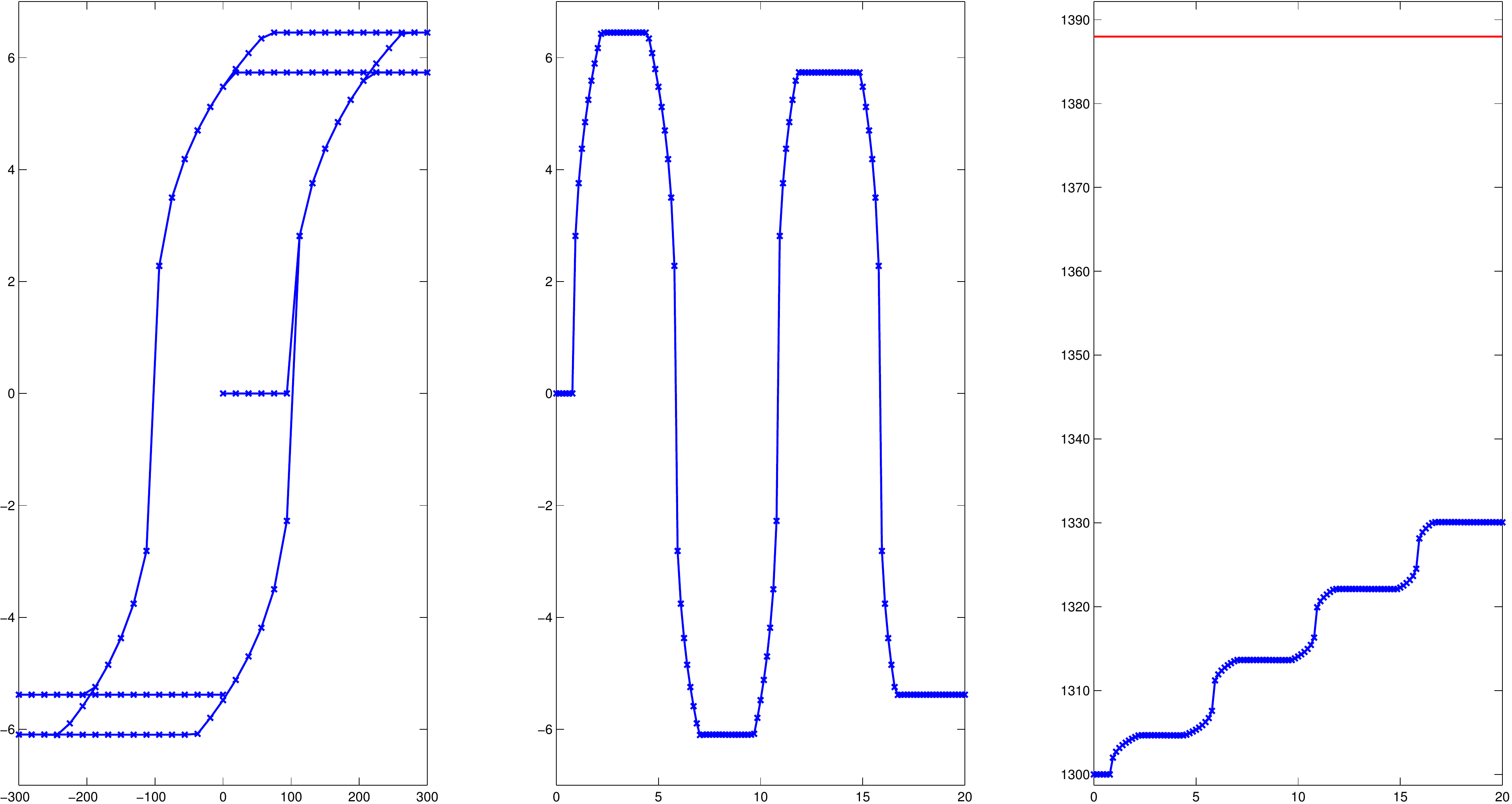}
\vspace{-0.7cm}
\caption{Average values of fields after two cycles of external forces: magnetization in x-direction versus external field (left), 
magnetization in x-direction versus time (middle), temperature versus time (right) never reaching the Curie temperature indicated by the red horizontal line.} 
\label{figure-two_loop}
\vspace{0.5cm}
\includegraphics[width=\textwidth,height=0.25\textheight]{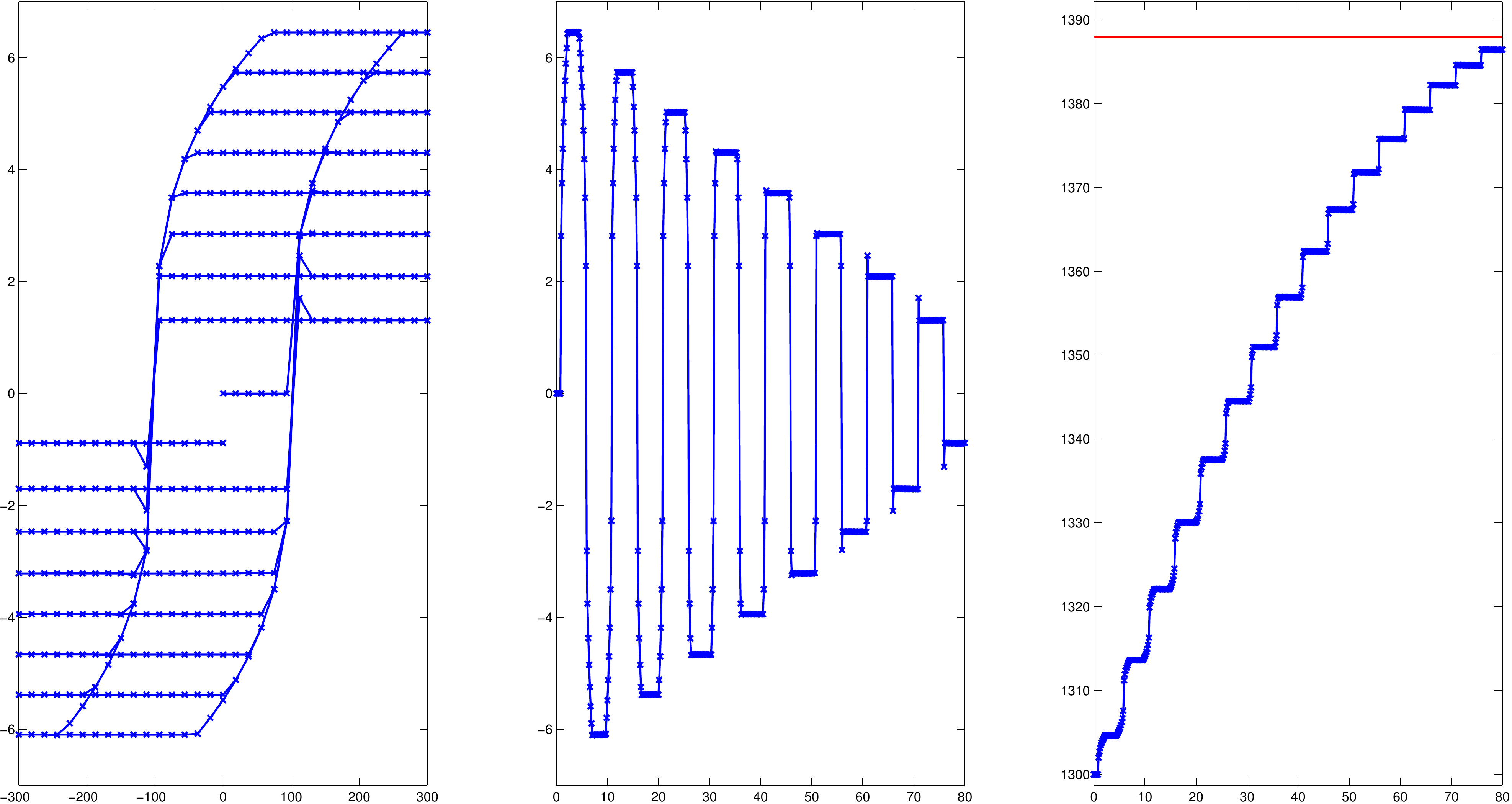}
\vspace{-0.7cm}
\caption{Average values of fields after eight cycles of external forces: magnetization in x-direction versus external field (left), 
magnetization in x-direction versus time (middle), temperature versus time (right)never reaching  the Curie temperature indicated by the red horizontal line.} 
\label{figure-eight_loops}
\end{figure}

\begin{figure}
\center
\includegraphics[width=\textwidth,height=0.25\textheight]{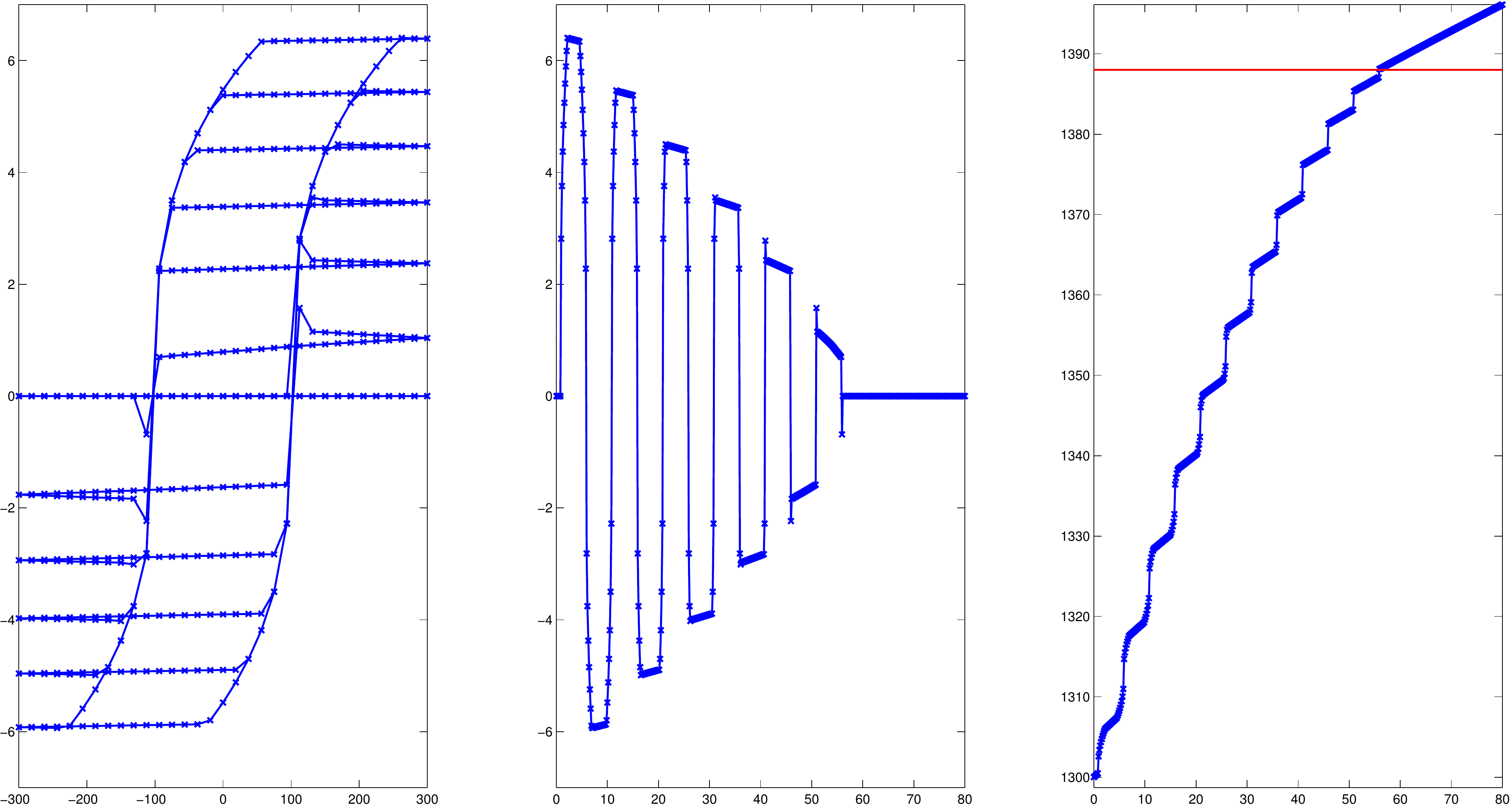}
\vspace{-0.7cm}
\caption{Average values of fields after eight cycles of external forces: magnetization in x-direction versus external field (left), 
magnetization in x-direction versus time (middle), temperature versus time (right) reaching and exceeding the Curie temperature indicated by the red horizontal line.} 
\label{figure-eight_loops_experiment2}
\end{figure}

As the result of the change of magnetic field inside the magnet, the magnet is heated and inside temperature increases with the boundary temperature $\theta_\mathrm{ext}$ held constant over time. An increase of the temperature decreases the measure support $p$, and amplitudes of magnetization become smaller over time. Figures \ref{figure-one_loop}-\ref{figure-eight_loops} describe average values of magnetization in $x$-direction and the temperature after one, two or eight cycles of external forces. With each cycle, the average temperature increases and approaches the Curie temperature. Since $\theta_\mathrm{ext} < \theta_c$, the temperature inside magnet never exceeds the Curie temperature and no paramagnetic effects are observed. A similar computation can be run with two modified physical parameters, $\theta_\mathrm{ext}=1500{\rm K}$, $b_0=0.1{\rm  W}/({\rm mK})$ . Then, the external temperature $\theta_\mathrm{ext} > \theta_c$  allows for heating up the magnet after the Curie temperature and a higher value of $b_0$ speeds up the heating process, see Figure \ref{figure-eight_loops_experiment2} for details.  It should be mentioned that choosing only $N=12$ Dirac atoms placed on 
``the middle sphere'' does not visibly change the shapes of Figures \ref{figure-one_loop}-\ref{figure-eight_loops}.

 The own MATLAB code is available as a package ``Thermo-magnetic solver'' at MATLAB Central and it can be downloaded for testing
at  \url{http://www.mathworks.com/matlabcentral/fileexchange/47878}.  It  utilizes the codes for an assembly of stifness and mass matrices described in \cite{rahman-valdman}. The assembly is vectorized and works very fast even for fine mesh triangulations. The inbuilt  MATLAB function $fmincon$ (it is a part of the Optimization Toolbox that must be available) was exploited for the minimization of \eqref{BalanceEqDis}. The function $fmincon$ was run with an automatic differentiation option, which is very time consuming even on coarse mesh triangulations. In order to speed up calculations of the magnetostatic potential $u_{m_{l,\tau}^k}$ from the Poisson problem \eqref{MaxwellDis}, an explicit inverse of the stiffness matrix was precomputed and stored for considered coarse mesh triangulations. Geometrical and material parameters can be adjusted for own testing in the functions $start.m$ and $start\_magnet.m$. \\

\section{Concluding remarks}\label{conclusion}
We tested computational performance of the model from \cite{benesova} on two-dimensional examples. In spite of a few simplifications (in particular, setting $\varkappa:=+\infty$), computational results are in qualitative agreement with physically observed phenomena. Interested readers are invited to perform their own numerical tests with a MATLAB code available on the web-page mentioned above. Adaptive approaches similar to the one in \cite{carstensen-roubicek,kruzik-prohl-1} could be used to allow for much finer discretizations of Young measure support and, as a consequence, for more accurate numerical approximations. Investigations of a convergence of the above scheme as well as verification of  discrete energy inequalities from \eqref{magneticbalance} and \eqref{thermalbalance} are left for our future work. 

\bigskip
\noindent
{\bf Acknowledgment:}  We thank anonymous referees for valuable comments and remarks which improved final exposition of our work.  We also  acknowledge the support by GA\v{C}R through projects 13-18652S, 16-34894L and by M\v{S}MT \v{C}R through project  7AMB16AT015.

\section{Appendix -- Young measures}
The Young
measures on a bounded  domain $\O\subset\R^n$ are weakly* measurable mappings
$x\mapsto\nu_x:\O\to \rca(\R^d)$ with values in probability measures;
 and the adjective ``weakly* measurable'' means that,
for any $v\in C_0(\R^d)$, the mapping
$\O\to\R:x\mapsto\A{\nu_x,v}=\int_{\R^d} v(\lambda)\nu_x(\d\lambda)$ is
measurable in the usual sense. Let us remind that, by the Riesz theorem,
$\rca(\R^d)$, normed by the total variation, is a Banach space which is
isometrically isomorphic with $C_0(\R^d)^*$, where $C_0(\R^d)$ stands for
the space of all continuous functions $\R^d\to\R$ vanishing at infinity.
Let us denote the set of all Young measures by $\mathscr{ Y}(\O;\R^d)$. It
is known that $\mathscr{Y}(\O;\R^d)$ is a convex subset of $L^\infty_{\rm
w}(\O;\rca(\R^d))\cong L^1(\O;C_0(\R^d))^*$, where the subscript ``w''
indicates the property ``weakly* measurable''.  A classical result
\cite{young} is that, for every sequence $\{y_k\}_{k\in\N}$
bounded in $L^\infty(\O;\R^d)$, there exists its subsequence (denoted by
the same indices for notational simplicity) and a Young measure
$\nu=\{\nu_x\}_{x\in\O}\in\mathscr{ Y}(\O;\R^d)$ such that
\be\label{jedna2}
\forall f\in C_0(\R^d):\ \ \ \ \lim_{k\to\infty}f\circ y_k=f_\nu\ \ \ \
\ \ \mbox{ weakly* in }L^\infty(\O)\ ,
\ee
where $[f\circ y_k](x)=f(y_k(x))$ and
\be
f_\nu(x)=\int_{\R^d}f(s)\nu_x(\d s)\ .
\ee
Let us denote by $\mathscr{Y}^\infty(\O;\R^d)$ the
set of all Young measures which are created by this way, i.e. by taking
all bounded sequences in $L^\infty(\O;\R^d)$. Note that (\ref{jedna2}) actually holds for
any $f:\R^d\to\R$ continuous.

A generalization of this result was formulated by
Schonbek \cite{schonbek} (cf. also \cite{roubicek0}): if
$1\le p<+\infty$: for every sequence
$\{y_k\}_{k\in\N}$ bounded in $L^p(\O;\R^{d})$ there exists its
subsequence (denoted by the same
indices) and a Young measure
$\nu=\{\nu_x\}_{x\in\O}\in\mathscr{Y}(\O;\R^d)$ such that
\be\label{young}
\forall f\in C_p(\R^d):\ \ \ \ \lim_{k\to\infty}f\circ y_k=f_\nu\
\ \ \ \ \ \mbox{ weakly in }L^1(\O)\ .\ee
We say that $\{y_k\}$ generates $\nu$ if (\ref{young}) holds.
Here for $p\ge 1$, we recall that   $C_p(\R^d)=\{f\in C(\R^d);\, f/(1+|\cdot|^p)\in C_0(\R^d)\}$.

  Let us denote by 
	$\mathscr{Y}^p(\O;\R^d)$ the set of all Young measures which are created by this
way, i.e. by taking all bounded sequences in $L^p(\O;\R^d)$.
It is well-known, however, that for any $\nu\in \mathscr{
Y}^p(\O;\R^d)$ there exists a special  generating sequence $\{y_k\}$ such that 
\eqref{young} holds even for $f\in C^p(\R^d)=\{y\in C(\R^d);\, |y|/(1+|\cdot|^p)\le C,\, C>0\}$.

\end{document}